\documentclass[amstex,12pt,amssymb]{article}
\usepackage{mathtext}
\usepackage{graphicx}
\usepackage{amsmath}
\usepackage{amssymb}
\usepackage{amsxtra}
\usepackage{latexsym}
\usepackage{ifthen}

\usepackage{indentfirst}

\usepackage{color}

\def\DD{P\ r\ o\ o\ f}

\makeatletter
\def\thmstyle{\it} 
\def\@begintheorem#1#2{\it \trivlist \item[\hskip
        \labelsep{\bf #1\ #2.}]\thmstyle}
\def\@opargbegintheorem#1#2#3{\it \trivlist \item[\hskip
        \labelsep{\bf #1\ #2\ (#3).}]\thmstyle}
\makeatother

\begin{document}

\renewcommand{\baselinestretch}{1.25}
\newcommand {\beq}{\begin{equation}}
\newcommand {\eeq}{\end{equation}}

\newtheorem{theorem}{Theorem}[section]
\newtheorem{lemma}[theorem]{Lemma}
\newtheorem{proposition}[theorem]{Proposition}
\newtheorem{corollary}[theorem]{Corollary}
\newtheorem{conjecture}[theorem]{Conjecture}
\newtheorem{definition}[theorem]{Definition}
\newtheorem{remark}[theorem]{Remark}
\newtheorem{claim}[theorem]{Claim}

\title{Large deviations in a population dynamics with catastrophes
 \author{\bf A. Logachov$^{1,2,3}$, O. Logachova$^{3,4}$,  A. Yambartsev$^{5}$
 }
}\maketitle

{\footnotesize
\noindent $^1$ Laboratory of Probability Theory and Mathematical Statistics, Sobolev Institute of Mathematrics,
Siberiian Branch of the RAS, Koptyuga str. 4, Novosibirsk 630090, RF\\
E-mail:\ \ avlogachov@mail.ru

\noindent $^2$ Novosibirsk State University,  Pirogova str. 1, Novosibirsk 630090, RF

\noindent $^3$ Novosibirsk State University of Economics and Management, Kamenskaya str. 56, Novosibirsk 630099, RF

\noindent $^4$ Siberian State University of Geo-systems and Technologies, Plakhotnogo str. 
10, 630108 Novosibirsk, RF\\
E-mail:\ \ omboldovskaya@mail.ru

\noindent $^5$ Institute of Mathematics and Statistics, University of S\~ao Paulo,
1010 Rua do Mat\~ao, CEP 05508--090, S\~ao Paulo SP, Brazil \\
E-mail:\ \ yambar@ime.usp.br}

\begin{abstract}
The large deviation principle on phase space is proved for a class of Markov processes
known as random population dynamics with catastrophes. In the paper we study the process 
which corresponds to the random population dynamics with linear growth and uniform catastrophes, 
where an eliminating portion of the population is chosen uniformly. 
The large deviation result provides an optimal trajectory of large fluctuation: 
it shows how the large fluctuations occur for this class of processes.   
\end{abstract}



\section{Introduction}

Stochastic models with catastrophes are studied since 70's and recieved a great attention of probability
community, 
see \cite{BGR1982} for the probably first systematic 
review about such processes, and see \cite{BRS2017} for short historical overview and more references. 
These models are used in analyzing a growth of a population subject to catastrophes due large-scale death
or emigrations of a population.
According \cite{BGR1982} the population dynamics considered here we will call population dynamics \textit{with linear growth and uniform catastrophes}, where the eliminating portion is chosen uniformly. 
Typically researchers are interested in  
extinction probability, the mean time to extinction, invariant measures, convergence to invariant measures for these processes. 


 In \cite{LLY1} for the population dynamic, $\xi(t)$ defined by (\ref{process0}, \ref{process}) in the following, with linear growth and uniform catastrophes we proved the local large deviation principle (LLDP): we established a rough 
logarithmic asymptotic for the probability of the scaling version $\xi_T(t), t\in[0,1]$, defined by 
(\ref{process.scaled}), 
to be in a small neighborhood of a continuous function. 
Here, based on the work \cite{LLY1} we established a large deviation principle (LDP) on the state space 
at the end of the interval of
observation of the process: we find the logarithmic asymptotic for the probability ${\bf P}(\xi_T(1)>x)$. 
Moreover, our proof also provides an optimal trajectory -- how such deviation 
occurs taking in account the evolution of the process.
As far as our understanding there are no other large deviations results for such processes. 

\vspace{0.5cm}

Throughout the paper we assume that all random elements are defined on a probability space
$(\Omega,\mathfrak{F},\bf{P})$.

\vspace{0.5cm}

We construct our continuous time process $\xi(t), t\in \mathbb{R}^+$, in two steps. 
First, consider the discrete time Markov chain $\eta(k)$, 
$k\in \mathbb{Z}^+$, $\mathbb{Z}^+=\{0\}\cup\mathbb{N}$,
with state space $\mathbb{Z}^+$ and transition probabilities 
\beq\label{process0}
{\bf P}(\eta(k+1)=j|\eta(k)=i)=
\left\{
\begin{array}{cl}
\frac{\lambda}{\lambda+\mu}, & \mbox{ if } j=i+1,\\
\frac{\mu}{i(\lambda+\mu)}, &\mbox{ if } 0\leq j<i, i\neq 0,\\
1, & \mbox{ if } j=1, i=0,
\end{array}
\right.
\eeq
where $\lambda$ and $\mu$ are positive constants. Let $\eta(0)=0$. Second, let $\nu(t)$, $t\in \mathbb{R}^+$, be 
the Poisson process with parameter $\mathbf{E}\nu(t)=\alpha t$, where $\alpha$ is positive parameter. 
Suppose that the process $\nu(\cdot)$ and the chain 
$\eta(\cdot)$ are independent. Finally,
\begin{equation}\label{process}
\xi(t):=\eta(\nu(t)), \ t\in \mathbb{R}^+.
\end{equation}

In order to establish the large deviation result we consider the scaled process
\begin{equation}\label{process.scaled}
\xi_T(t):=\frac{\xi(Tt)}{T}, \ t\in [0,1], \mbox{ as }T\to\infty.
\end{equation}
where $T$ is an increasing parameter, $T\to\infty$. We are interested in LDP for the family 
of random variables $\xi_T(1)$.

In the proof of LDP we use the standard implication (see, for example \cite[Lemma 4.1.23]{DZ}):
$$
\text{LLDP and ET} \Rightarrow \text{LDP},
$$
where LLDP is Local Large Deviation Principle, and ET stands for Exponential tightness.

In the next section we recall definitions and formulate main result. In Section~\ref{sec2} we prove the main result, 
Theorem~\ref{t2.2}. Auxiliary results are proved in Section~\ref{sec3}.

\section{Definitions and main results}

Recall the definitions we need.
\begin{definition}
A family of random variables $\xi_T(1)$ satisfies
LLDP  in  $\mathbb{R}$ with the rate function
$I = I(x):\ \mathbb{R} \rightarrow [0,\infty]$ and the normalizing function
$\psi(T):\ \lim\limits_{T\rightarrow\infty}\psi(T) = \infty$,
 if the following equality holds for any  $x \in  \mathbb{R} $
$$
\lim_{\varepsilon\rightarrow 0}\limsup_{T\rightarrow \infty}\frac{1}{\psi(T)}
\ln\mathbf{P}(\xi_T(1)\in U_\varepsilon(x))
=\lim_{\varepsilon\rightarrow 0}\liminf_{T\rightarrow \infty}\frac{1}{\psi(T)}
\ln\mathbf{P}(\xi_T(1)\in U_\varepsilon(x))=-I(x),
$$
where
$
U_\varepsilon(x):=\{y\in \mathbb{R}: \ |x-y|<\varepsilon\}.
$
\end{definition}

\begin{definition} \label{d1.2}
  A family of random variables $\xi_T(1)$ is
exponentially tight  on  $\mathbb{R}$  if, for any $c < \infty$, there exists a compact set $K_c \in \mathbb{R}$ 
such that
$$
 \limsup_{T \rightarrow \infty} \frac{1}{\psi(T)} \ln \mathbf{P}(\xi_T(1) \not\in K_c  ) <-c.
$$
\end{definition}
We denote the closure and interior of the set $B$ by  $[B]$ and $(B)$, respectively.

\begin{definition} \label{d1.3}
  A family of random variables $\xi_T(1)$ satisfies
LDP  on  $\mathbb{R}$ with a rate function
$I = I(f)\,:\, \mathbb{R} \rightarrow [0,\infty]$ and the normalizing function
$\psi(T):\ \lim\limits_{T\rightarrow\infty}\psi(T) = \infty$,
if for any $c \geq 0 $ \ the set $\{ x \in \mathbb{R}:\ I(x) \leq c \}$  is a compact set 
 and, for any set $B \in \mathfrak{B}(\mathbb{R})$ the following inequalities hold:
$$
\begin{aligned}
 \limsup_{T \rightarrow \infty} \frac{1}{\psi(T)} \ln \mathbf{P}(\xi_T(1) \in B ) &\leq - I([B]), \\
\liminf_{T \rightarrow \infty} \frac{1}{\psi(T)} \ln \mathbf{P}(\xi_T(1) \in B ) &\geq -I((B)),
\end{aligned}
$$
where $\mathfrak{B}(\mathbb{R})$ is the Borel $\sigma$-algebra on $\mathbb{R}$, $I(B) = \inf\limits_{x \in B} I(x)$,
$I(\emptyset) = \infty$.
\end{definition}

Further we will use the following notations: $\overline{B}$ is a complement of the set $B$;
$\mathbf{I}(B)$ is the indicator function of the set $B$;
$[a]$ is the integer part of the number $a$.

We recall Low of Large Numbers (LLN), it was proved in \cite{LLY1}.

\begin{theorem} (LLN) \label{t2.1}  For any $\varepsilon>0$ the following equality holds true
$$
\mathbf{P}\Bigl(\lim\limits_{T\rightarrow\infty}\sup\limits_{t\in[0,1]}\xi_T(t)>\varepsilon\Bigr)=0.
$$
\end{theorem}
LDP is the main theorem in the paper.

\begin{theorem} (LDP) \label{t2.2}  The family of random variables $\xi_T(1)$ satisfies LDP with 
normalized function 
$\psi(T)=T$ and with rate function
$$
I(x)= \left\{
           \begin{array}{ll} 
                              \ \infty,  & \text{ if } x\in (-\infty,0),\\
                              x\ln\left(\frac{\lambda+\mu}{\lambda}\right), & \text{ if } x\in [0,\alpha),\\
                              x\ln\left(\frac{x(\lambda+\mu)}{\alpha\lambda}\right)-x+\alpha, & \text{ if } x\in [\alpha,\infty).\\
                              \end{array}
                              \right.
$$
\end{theorem}

\begin{figure}[ht!]
\begin{center}
\includegraphics[width=0.6\linewidth]{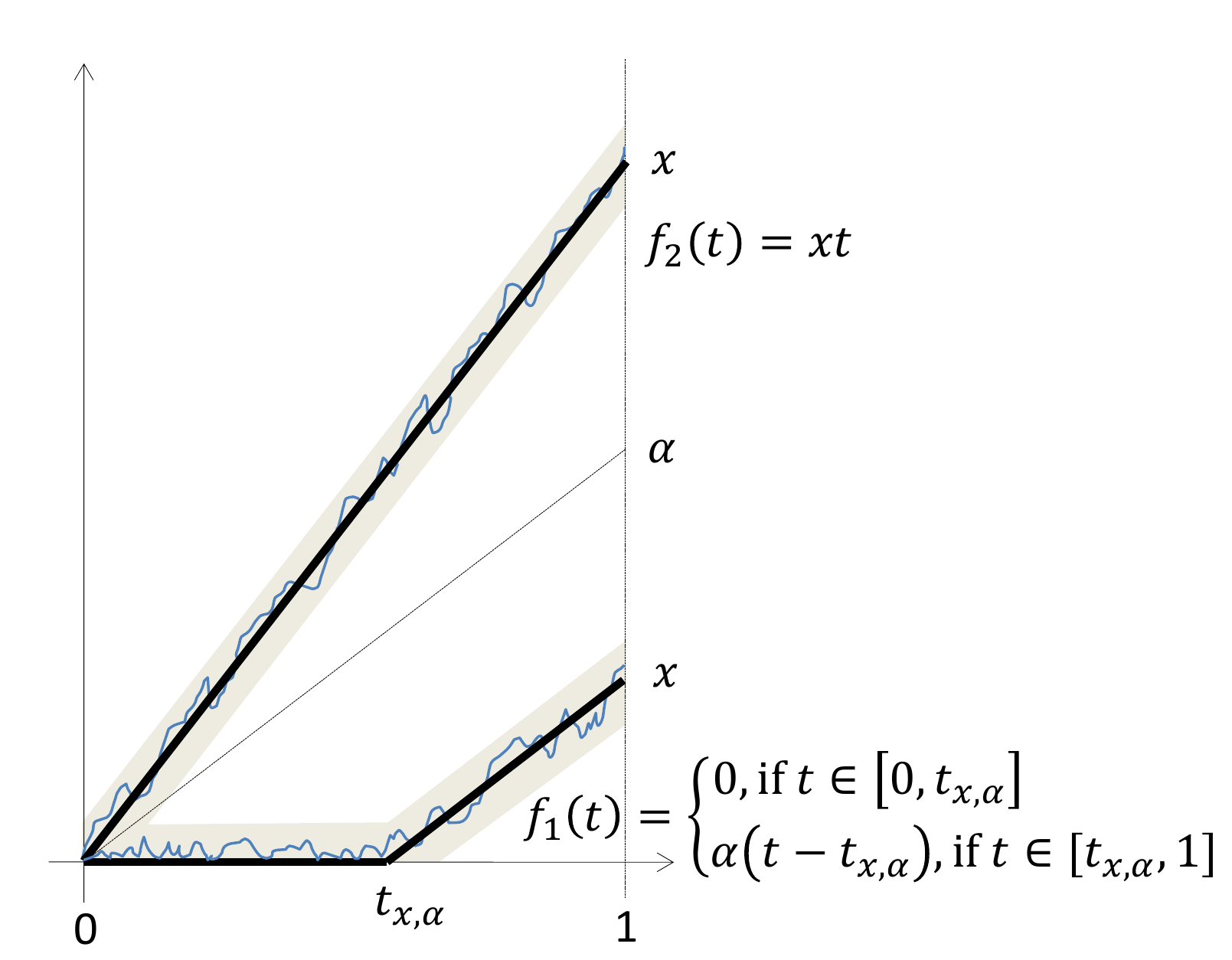}
\end{center}
\caption { The ``most probable'' trajectories which provide large fluctuations. 
If $x < \alpha$ then the large deviation occurs according the function $f_1$. If $x\ge \alpha$, then the large deviation trajectory is in the neighborhood of the straight line $f_2$.}%
\label{fig1}
\end{figure}

The proof of Theorem~\ref{t2.2} provides the ``most probable'' trajectories of large deviations $\xi_T(1) > x$. If $x < \alpha$ then
there exists the moment $t_{x,\alpha} = 1- \frac{x}{\alpha} \in (0,1)$ such that the process $\xi_T(\cdot)$ stays 
near the zero up to the time $t_{x,\alpha}$ and after that $\xi_T(t), t\ge t_{x,\alpha}$ increases according  
the straight line
which starts at point $(t_{x,\alpha}, 0)$ and grows up to the point $(1,x)$ with the slope $\alpha$, 
see the function $f_1$ on Figure~\ref{fig1}. If $x \ge \alpha$ then 
the process grows together with the straight line starting from origin up to the point $(1,x)$, i.e. its slope is $x$, the function $f_2$ on Figure~\ref{fig1}.

\section{Proof of Theorem~\ref{t2.2}}\label{sec2}

\noindent
Random process 
$\xi(t)$ we represent in the following way
\beq \label{4.1}
\xi(t)=\nu_1(t)-\sum\limits_{k=0}^{\nu_2(t)}\zeta_k(\xi(\tau_k-)),
\eeq
where $\nu_1(t)$, $\nu_2(t)$ are independent Poisson point processes
with parameters $$\mathbf{E}\nu_1(t)=\dfrac{\alpha\lambda}{\lambda+\mu} t, \ \
\mathbf{E}\nu_2(t)=\dfrac{\alpha\mu}{\lambda+\mu} t;$$
$0= \tau_0 < \tau_1<\dots<\tau_k\dots$ -- jump moments of the process $\nu_2(t)$;
random variables $\zeta_k(m)$, $k\in \mathbb{Z}^+$, $m\in \mathbb{Z}^+$ are mutually independent
and do not depend on $\nu_1(t)$ and $\nu_2(t)$; $\zeta_0(m)=0$, for all $m\in \mathbb{Z}^+$; for fixed
$k,m \in \mathbb{N}$ the distribution of $\zeta_k(m)$ is given by
$$
\mathbf{P}(\zeta_k(m)=r)= \frac 1 m, \ 1\leq r\leq m,
$$
and $\zeta_k(0)=-1$, for all $k\in \mathbb{N}$.
Using representation (\ref{4.1}), we have
$$
\xi_T(t)=\frac{\nu_1(Tt)}{T}-\frac{1}{T}\sum\limits_{k=0}^{\nu_2(Tt)}\zeta_k(\xi(\tau_k-)):=\xi_{T}^+(t)-\xi_{T}^-(t).
$$
From Theorem~\ref{t2.1} follows that Theorem~\ref{t2.2} holds true for $x=0$. 
Prove the theorem for $x>0$.
Let us estimate from above $\mathbf{P}(\xi_T(1)\geq x)$, $x>0$.
For any $\delta>0$ and for any $n\in \mathbb{N}$ we obtain
\beq \label{4.2}
\begin{aligned}
&\mathbf{P}(\xi_T(1)\geq x)=\mathbf{P}\Bigr(\nu_1(T)+\sum\limits_{k=0}^{\nu_2(T)}\zeta_k(\xi(\tau_k-))\geq Tx\Bigr)
\\
&\leq \sum\limits_{k=1}^{n-1}
\mathbf{P}\Bigl(\nu_1(T)+\sum\limits_{k=0}^{\nu_2(T)}\zeta_k(\xi(\tau_k-))\geq Tx
,\sup\limits_{t\in\left[\frac{k-1}{n},\frac{k}{n}\right]}\xi_T(t)\leq\delta,
\sup\limits_{t\in\left(\frac{k}{n},1\right]}\xi_T(t)>\delta\bigr)
\\
&
+\mathbf{P}\Bigl(\nu_1(T)+\sum\limits_{k=0}^{\nu_2(T)}\zeta_k(\xi(\tau_k-))\geq Tx
,\sup\limits_{t\in\left[\frac{n-1}{n},1\right]}\xi_T(t)\leq\delta\Bigr):=\mathbf{P}_1+\mathbf{P}_2.
\end{aligned}
\eeq
Estimate from above
$\mathbf{P}_1$.
\beq \label{4.3}
\mathbf{P}_1\leq (n-1)\max\limits_{1\leq k \leq n-1}
\mathbf{P}\Bigl(A,
\inf\limits_{t\in\left[\frac{k-1}{n},\frac{k}{n}\right]}\xi_T(t)\leq\delta,
\inf\limits_{t\in\left(\frac{k}{n},1\right]}\xi_T(t)>\delta\Bigr):=(n-1)\max\limits_{1\leq k \leq n-1}
\mathbf{P}_{1k},
\eeq
where 
$$
A:=\Bigr\{\omega:\nu_1(T)+\sum\limits_{k=0}^{\nu_2(T)}\zeta_k(\xi(\tau_k-))\geq Tx\Bigl\}.
$$
Estimate from above $\mathbf{P}_{1k}$. For any $c>0$ we have
\beq \label{4.4}
\begin{aligned}
\mathbf{P}_{1k} & \leq
\mathbf{P}\Bigr(\nu_1(T)-\nu_1\Bigr(T\frac{k-1}{n}\Bigl)\geq T(x-\delta),B_c\Bigl)
\\
&+\mathbf{P}\Bigr(A,\inf\limits_{t\in \left[\frac{k-1}{n},\frac{k}{n}\right]}\xi_T(t)\leq\delta,
\inf\limits_{t\in\left(\frac{k}{n},1\right]}\xi_T(t)>\delta,\overline{B_c}\Bigl):=\mathbf{P}_{1k}^1+\mathbf{P}_{1k}^2,
\end{aligned}
\eeq
where
$$
B_c:=\Bigl\{\omega:\nu_2(T)-\nu_2\Bigl(T\frac{k}{n}\Bigr)\leq cT\Bigr\}.
$$
Since $\nu_1(\cdot)$ and $\nu_2(\cdot)$ are independent, then
\beq \label{4.5}
\mathbf{P}_{1k}^1=
\mathbf{P}\Bigl(\nu_1(T)-\nu_1\Bigl(T\frac{k-1}{n}\Bigr)\geq T(x-\delta)\Bigr)\mathbf{P}(B_c).
\eeq
Estimate from above $\mathbf{P}_{1k}^2$. For any $a>0$ we obtain
\beq \label{4.6}
\begin{aligned}
& \mathbf{P}_{1k}^2 \leq\mathbf{P}\Bigl(A,\inf\limits_{t\in\left[\frac{k-1}{n},\frac{k}{n}\right]}\xi_T(t)\leq\delta,
\inf\limits_{t\in\left(\frac{k}{n},1\right]}\xi_T(t)>\delta,\overline{B_c},
\sum\limits_{l=1}^{[cT]}\zeta_k(\xi(\tau_{k_l}-))>aT\Bigr)
\\
& \qquad +\mathbf{P}\Bigl(\inf\limits_{t\in\left(\frac{k}{n},1\right]}\xi_T(t)>\delta,\overline{B_c},
\sum\limits_{l=1}^{[cT]}\zeta_k(\xi(\tau_{k_l}-))\leq aT\Bigr)
\\
&\leq \mathbf{P}\Bigl(\nu_1(T)-\nu_1\Bigl(T\frac{k-1}{n}\Bigr)\geq T(x+a-\delta)\Bigr)
\\
& \qquad +\mathbf{P}\Bigl(\inf\limits_{t\in\left(\frac{k}{n},1\right]}\xi_T(t)>\delta,\overline{B_c},
\sum\limits_{l=1}^{[cT]}\zeta_{k_l}(\xi(\tau_{k_l}-))\leq aT\Bigr)
\\
& \leq \mathbf{P}\Bigl(\nu_1(T)-\nu_1\Bigl(T\frac{k-1}{n}\Bigr)\geq T(x+a-\delta)\Bigr)
+\mathbf{P}\Bigl(\bigcap\limits_{l=1}^{[cT]}C_l,\overline{B_c},
\sum\limits_{l=1}^{[cT]}\zeta_{k_l}(\xi(\tau_{k_l}-))\leq aT\Bigr),
\end{aligned}
\eeq
where
$\tau_{k_1},\dots,\tau_{k_{[cT]}}$ are the first $[cT]$ jump times 
of the process $\nu_2(Tt)$ for $t\in \big[\frac{k}{n},1\big]$, and
\begin{equation}\label{Cl}
C_l:=\{\omega:\xi_T(\tau_{k_l}-)>\delta\}.
\end{equation}
For $T$ sufficiently large and $x+a-\delta \ge 1$ we obtain 
\beq \label{4.7}
\begin{aligned}
\mathbf{P} & \Bigl(\nu_1(T)-\nu_1\Bigl(T\frac{k-1}{n}\Bigr)\geq T(x+a-\delta)\Bigr)
=
 e^{-\frac{\alpha\lambda(n-k+1)}{n(\lambda+\mu)}T}
\!\!\!\!\!\!
\sum\limits_{r=[T(x+a-\delta)]+1}^\infty
\!\!\!\!\!\! \frac{(\alpha\lambda (n-k+1) T)^r}{r!n^r(\lambda+\mu)^r}
\\
&\leq e^{-\frac{\alpha\lambda(n-k+1)}{n(\lambda+\mu)}T}
\Bigl(\frac{\alpha\lambda (n-k+1) T}{[T(x+a/2-\delta)]n(\lambda+\mu)}\Bigr)^{[Ta]/2}
\!\!\!\!\!\!
\sum\limits_{r=[T(x+a/2-\delta)]+1}^\infty
\!\!\!\!\!\!
\frac{(\alpha\lambda (n-k+1) T)^r}{r!n^r(\lambda+\mu)^r}
\\
&\leq e^{-\frac{\alpha\lambda(n-k+1)}{n(\lambda+\mu)}T}
\Bigl(\frac{2\alpha\lambda}{[(x+a/2-\delta)](\lambda+\mu)}\Bigr)^{[Ta]/2}
\!\!\!\!\!\!
\sum\limits_{r=[T(x-\delta)]+1}^\infty
\!\!\!\!\!\!
\frac{(\alpha\lambda (n-k+1) T)^r}{r!n^r(\lambda+\mu)^r}.
\end{aligned}
\eeq
Since there exists $a^*(x,\alpha,\delta,\lambda,\mu)$ such that for all $a\geq a^*$ the inequality
$$
\bigg(\frac{2\alpha\lambda}{[(x+a/2-\delta)](\lambda+\mu)}\bigg)^{[Ta]/2}\leq e^{-\frac{\alpha\mu}{(\lambda+\mu)}T}\leq  \mathbf{P}(B_c)
$$
holds true, then from inequality (\ref{4.7}) it follows that for $a>a^*$
\beq \label{4.8}
\mathbf{P}\bigg(\nu_1(T)-\nu_1\bigg(T\frac{k-1}{n}\bigg)\geq T(x+a-\delta)\bigg)\leq\mathbf{P}_{1k}^1.
\eeq
From Lemma~\ref{l5.5} it follows that for any $a>0$ and for $T$ sufficiently large
\beq \label{4.9}
\mathbf{P}\bigg(\bigcap\limits_{l=1}^{[cT]}C_l, \overline{B_c},
\sum\limits_{l=1}^{[cT]}\zeta_{k_l}(\xi(\tau_{k_l}-))\leq aT\bigg)\leq
\bigg(\frac{1}{[\delta T]}\bigg)^{[cT]}\exp(aT)\leq \mathbf{P}_{1k}^1.
\eeq
Choose $a>a^*$. Using (\ref{4.4})--(\ref{4.7}), (\ref{4.8}), (\ref{4.9}), we obtain 
for sufficiently large $T$
\beq \label{4.10}
\mathbf{P}_{1k}\leq 3\mathbf{P}_{1k}^1
=3\mathbf{P}\bigg(\nu_1(T)-\nu_1\bigg(T\frac{k-1}{n}\bigg)\geq T(x-\delta)\bigg)\mathbf{P}(B_c).
\eeq
From (\ref{4.3}), (\ref{4.10}) it follows that for all 
$n\in \mathbb{N}$, $\delta>0$, $c>0$
\beq \label{4.11}
\mathbf{P}_1\leq 3(n-1)\max\limits_{1\leq k \leq n-1}\mathbf{P}\bigg(\nu_1(T)-\nu_1\bigg(T\frac{k-1}{n}\bigg)\geq T(x-\delta)\bigg)\mathbf{P}(B_c).
\eeq
Estimate from above  $\mathbf{P}_2$. 
The following inequalities holds true for $n\geq \exp\big(\frac{\alpha}{x-\delta}\big)$
$$
\begin{aligned}
\mathbf{P}_2&=\mathbf{P}\Bigl(\nu_1(T)+\sum\limits_{k=0}^{\nu_2(T)}\zeta_k(\xi(\tau_k-))\geq Tx
,\inf\limits_{t\in\left[\frac{n-1}{n},1\right]}\xi_T(t)\leq\delta\Bigr)
\\
&\leq\mathbf{P}\Bigl(\nu_1(T)-\nu_1\Bigl(T\frac{n-1}{n}\Bigr)\geq T(x-\delta)\Bigr)=
e^{-\frac{\alpha\lambda}{n(\lambda+\mu)}T}\sum\limits_{k=[T(x-\delta)]+1}^\infty
\frac{(\alpha\lambda T)^k}{k!(n(\lambda+\mu))^k}
\\
&\leq e^{-\frac{\alpha\lambda}{n(\lambda+\mu)}T}\bigg(\frac{1}{n}\bigg)^{T(x-\delta)}\sum\limits_{k=[T(x-\delta)]+1}^\infty
\frac{(\alpha\lambda T)^k}{k!(\lambda+\mu)^k}
\leq e^{-\alpha T}
\sum\limits_{k=[T(x-\delta)]+1}^\infty
\frac{(\alpha\lambda T)^k}{k!(\lambda+\mu)^k}
\\
&= e^{-\frac{\alpha\lambda}{(\lambda+\mu)}T}e^{-\frac{\alpha\mu}{(\lambda+\mu)}T}\sum\limits_{k=[T(x-\delta)]+1}^\infty
\frac{(\alpha\lambda T)^k}{k!(\lambda+\mu)^k}\leq \mathbf{P}_{11}^1.
\end{aligned}
$$
Thus, from (\ref{4.2}), (\ref{4.11}) it follows that for $n\geq \exp\big(\frac{\alpha}{x-\delta}\big)$
\beq \label{4.12}
\mathbf{P}(\xi_T(1)\geq x)\leq 4(n-1) \max\limits_{1\leq k \leq n-1}
\mathbf{P}\Bigl(\nu_1(T)-\nu_1\Bigl(T\frac{k-1}{n}\Bigr)\geq T(x-\delta)\Bigr)\mathbf{P}(B_c).
\eeq
Using Lemma~\ref{l5.3}, Lemma~\ref{l5.4}, Definition~\ref{d1.3} and inequality (\ref{4.12}) we obtain for 
$n\geq \exp\big(\frac{\alpha}{x-\delta}\big)$, $\delta>0$, $c>0$
$$
\begin{aligned}
\limsup\limits_{T \rightarrow \infty}\frac{1}{T} \ln \mathbf{P}(\xi_T(1)\geq x)
\leq &
\max\limits_{1\leq k \leq n-1}\sup\limits_{y\geq x}\bigg(-(y-\delta)\ln\left(\frac{(y-\delta)(\lambda+\mu)n}{\alpha\lambda(n-k+1)}\right)
+(y-\delta)
\\
& {} -\frac{\alpha\lambda(n-k+1)}{(\lambda+\mu)n}-\frac{\alpha\mu(n-k)}{(\lambda+\mu)n}
 +\frac{\alpha\mu(n-k)c}{(\lambda+\mu)n}- c\ln c\bigg).
\end{aligned}
$$
When $\delta\rightarrow 0$, $c\rightarrow 0$ we obtain
\small
$$
\limsup\limits_{T \rightarrow \infty}\frac{1}{T} \ln \mathbf{P}(\xi_T(1)\geq x)
\leq
\max\limits_{1\leq k \leq n-1}\sup\limits_{y\geq x}\bigg(-y\ln\left(\frac{y(\lambda+\mu)n}{\alpha\lambda(n-k+1)}\right)
+y-\frac{\alpha\lambda(n-k+1)}{(\lambda+\mu)n}-\frac{\alpha\mu(n-k)}{(\lambda+\mu)n}\bigg).
$$
\normalsize
And when $n\rightarrow\infty$ we obtain
$$
\begin{aligned}
\limsup\limits_{T \rightarrow \infty}\frac{1}{T} \ln \mathbf{P}(\xi_T(1)\geq x)
&\leq
\sup\limits_{z\in(0,1)}\sup\limits_{y\geq x}\bigg(-y\ln\left(\frac{y(\lambda+\mu)}{\alpha\lambda z}\right)
+y-\frac{\alpha\lambda z}{(\lambda+\mu)}-\frac{\alpha\mu z}{(\lambda+\mu)}\bigg)
\\
&=\sup\limits_{z\in(0,1)}\sup\limits_{y\geq x}\bigg(-y\ln\left(\frac{y(\lambda+\mu)}{\alpha\lambda z}\right)
+y - \alpha z \bigg).
\end{aligned}
$$
Finding the maximum of the function
$$
f(y,z)=-y\ln\left(\frac{y(\lambda+\mu)}{\alpha\lambda z}\right)
+y - \alpha z
$$
in the domain $y\geq x$, $z\in[0,1]$ we obtain
\beq \label{4.13}
\limsup\limits_{T \rightarrow \infty}\frac{1}{T} \ln \mathbf{P}(\xi_T(1)\geq x)
\leq - I(x).
\eeq

Estimate from below $\mathbf{P}(\xi_T(1)>x)$. 
Since the processes $\nu_1(\cdot)$ and $\nu_2(\cdot)$ are independent, each with independent increments,
then 
for all $\varepsilon>0$, $z\in(0,1)$
$$
\begin{aligned}
&\mathbf{P}(\xi_T(1)>x)=\mathbf{P}\Big(\nu_1(T)+\sum\limits_{k=0}^{\nu_2(T)}\zeta_k(\xi(\tau_k-))\geq Tx\Big)
\\
&\geq \mathbf{P}\Big(\xi_T(1-z)<\varepsilon,\nu_1(T)-\nu_1(T(1-z))>Tx,\nu_2(T)-\nu_2(T(1-z))=0\Big)
\\
&= \mathbf{P}\big(\xi_T(1-z)<\varepsilon\big)\mathbf{P}\big(\nu_1(T)-\nu_1(T(1-z))>Tx\big)
\mathbf{P}\big(\nu_2(T)-\nu_2(T(1-z))=0\big).
\end{aligned}
$$
From Theorem~\ref{t2.1} and Lemma~\ref{l5.3} it follows that for any $z\in(0,1)$
$$
\liminf\limits_{T \rightarrow \infty}\frac{1}{T} \ln \mathbf{P}(\xi_T(1)> x)
\geq \sup\limits_{y\geq x}\bigg(-y\ln\left(\frac{y(\lambda+\mu)}{\alpha\lambda z}\right)
+y-\frac{\alpha\lambda z}{(\lambda+\mu)}\bigg)-\frac{\alpha\mu z}{(\lambda+\mu)}.
$$
Since it holds for any $z\in(0,1)$, then
\beq \label{4.14}
\liminf\limits_{T \rightarrow \infty}\frac{1}{T} \ln \mathbf{P}(\xi_T(1)> x)
\geq \sup\limits_{z\in(0,1)}\sup\limits_{y\geq x}\bigg(-y\ln\left(\frac{y(\lambda+\mu)}{\alpha\lambda z}\right)
+y - \alpha z \bigg)=-I(x).
\eeq

Let us prove the LLDP for the family of the random variables $\xi_T(1)$.
Applying (\ref{4.13}), we obtain
$$
\begin{aligned}
\lim\limits_{\varepsilon \rightarrow 0}\limsup\limits_{T \rightarrow \infty}\frac{1}{T} \ln \mathbf{P}\big(\xi_T(1)\in (x-\varepsilon,x+\varepsilon) \big)
&\leq\lim\limits_{\varepsilon \rightarrow 0}\limsup\limits_{T \rightarrow \infty}\frac{1}{T} \ln \mathbf{P}(\xi_T(1)\geq x-\varepsilon)
\\
&\leq-\lim\limits_{\varepsilon \rightarrow 0}I(x-\varepsilon)=-I(x).
\end{aligned}
$$
Using (\ref{4.13}), (\ref{4.14}) and the fact that the function $I(x)$ is continuously differentiable
function for $x>0$, we obtain
$$
\begin{aligned}
&\lim\limits_{\varepsilon \rightarrow 0}\liminf\limits_{T \rightarrow \infty}\frac{1}{T} \ln \mathbf{P}(\xi_T(1)\in (x-\varepsilon,x+\varepsilon))
\\
&=\lim\limits_{\varepsilon \rightarrow 0}\liminf\limits_{T \rightarrow \infty}\frac{1}{T} 
\ln \big(\mathbf{P}(\xi_T(1)> x-\varepsilon)-\mathbf{P}(\xi_T(1) \geq x+\varepsilon)\big)
\\
&=\lim\limits_{\varepsilon \rightarrow 0}\liminf\limits_{T \rightarrow \infty}\frac{1}{T}
\ln\mathbf{P}(\xi_T(1)> x-\varepsilon) + 
\lim\limits_{\varepsilon \rightarrow 0}\liminf\limits_{T \rightarrow \infty}\frac{1}{T}
\ln\bigg(1-\frac{\mathbf{P}(\xi_T(1) \geq x+\varepsilon)}{\mathbf{P}(\xi_T(1) > x-\varepsilon)}\bigg)
\\
&\geq \lim\limits_{\varepsilon \rightarrow 0}\liminf\limits_{T \rightarrow \infty}\frac{1}{T}
\ln\mathbf{P}(\xi_T(1)> x-\varepsilon)+
\lim\limits_{\varepsilon \rightarrow 0}\liminf\limits_{T \rightarrow \infty}\frac{1}{T}
\ln\bigg(1-e^{-T(I(x+\varepsilon)-I(x-\varepsilon)+o(1))}\bigg)
\\
&=-I(x)+\lim\limits_{\varepsilon \rightarrow 0}\liminf\limits_{T \rightarrow \infty}\frac{1}{T}
\ln\bigg(1-e^{-T(2\varepsilon I'(\tilde{x})+o(1))}\bigg)=-I(x),
\end{aligned}
$$
where $\tilde{x}\in(x-\varepsilon,x+\varepsilon)$
 is the point where $I'(\tilde{x})=(I(x+\varepsilon)-I(x-\varepsilon))/2\varepsilon$.
 Thus LLDP is proved.
 
Exponential tightness follows from (\ref{4.13}) and the fact that
for any $c\geq 0$ the set $I(x)\leq c$ is compact. $\Box$

\section{Auxiliary results}\label{sec3}

Here we will prove several auxiliary lemmas.

\begin{lemma} \label{l5.2}
For any $a\in \mathbb{R}$ the following inequality holds true
$$
\mathbf{P}\bigg(\sum\limits_{l=1}^{[cT]}\zeta_{k_l}([\delta T])\leq aT\bigg)\leq
\bigg(\frac{1}{[\delta T]}\bigg)^{[cT]}\exp(aT).
$$
\end{lemma}

\noindent
\DD. 
Since 
$\zeta_{k_l}([\delta T])$, $1\leq l \leq [cT]$ are i.i.d., using Chebyshev inequality we obtain
$$
\begin{aligned}
\mathbf{P}\bigg(\sum\limits_{l=1}^{[cT]}\zeta_{k_l}([\delta T])\leq aT\bigg)&=
\mathbf{P}\bigg(\exp\bigg(-\sum\limits_{l=1}^{[cT]}\zeta_{k_l}([\delta T])\bigg)\geq \exp(-aT)\bigg)
\\
&\leq\frac{\mathbf{E}\exp\bigg(-\sum\limits_{l=1}^{[cT]}\zeta_{k_l}([\delta T])\bigg)}{\exp(-aT)}=
\frac{\bigg(\mathbf{E}\exp(-\zeta_{k_1}([\delta T]))\bigg)^{[cT]}}{\exp(-aT)}
\\
&=\frac{\bigg(\frac{1}{[\delta T]}\sum\limits_{r=1}^{[\delta T]}\exp(-r)\bigg)^{[cT]}}{\exp(-aT)}
\leq \bigg(\frac{1}{[\delta T]}\bigg)^{[cT]}\exp(aT).
\end{aligned}
$$
$\Box$

\begin{lemma} \label{l5.3}
The family of random variables $\frac{1}{T}\big(\nu_1(T)-\nu_1(T\Delta)\big)$, $\Delta\in[0,1)$ satisfies LDP
with normalized function 
$\psi(T)=T$ and rate function
$$
I_1(x)= \left\{
           \begin{array}{ll}
           \infty,  & \text{ if }\ x\in (-\infty,0),\\
           x\ln\left(\frac{x(\lambda+\mu)}{\alpha\lambda(1-\Delta)}\right)-x+
           \frac{\alpha\lambda(1-\Delta)}{\lambda+\mu}, & \text{ if }\ x\in [0,\infty).
                              \end{array}
                              \right.
$$
\end{lemma}

\noindent
\DD. Random variable $\frac{1}{T}\big(\nu_1(T)-\nu_1(T\Delta)\big)$
 can be represented as a sum of independent random variables which have the same 
 distribution as $\nu_1(1-\Delta)$. Thus from \cite[Theorem 2.2.3]{DZ}
it is enough to show that the Legendre transform of exponential moment of random variable 
$\nu_1(1-\Delta)$ has the following form
$$
\Lambda(x)=\sup\limits_{y\in \mathbb{R}} \big( xy-\ln\mathbf{E}e^{y\nu_1(1-\Delta)} \big)=
x\ln\left(\frac{x(\lambda+\mu)}{\alpha\lambda(1-\Delta)}\right)-x+\frac{\alpha\lambda(1-\Delta)}{\lambda+\mu}, \ x\geq 0.
$$
Since
$$
\mathbf{E}e^{y\nu_1(1-\Delta)}=\exp\bigg\{\frac{\alpha\lambda(1-\Delta)}{\lambda+\mu} e^{y}-\frac{\alpha\lambda(1-\Delta)}{\lambda+\mu}\bigg\},
$$
then the differential calculus finishes the proof. $\Box$

\begin{lemma} \label{l5.4} For any $c\in[0,1)$, $\Delta\in [0,1]$ the following inequality holds true
\beq \label{5.4}
\mathbf{P} \big( \nu_2(T)-\nu_2(\Delta T)\leq cT \big)\leq \exp\Bigl\{-\frac{\alpha\mu(1-\Delta)}{\lambda+\mu}T
 +\frac{\alpha\mu(1-\Delta)c}{\lambda+\mu}T-T c\ln c \Bigr\}.
\eeq
\end{lemma}

\noindent
\DD. For any  $r>0$ by Chebyshev inequality we have
$$
\begin{aligned}
& \mathbf{P}\big( \nu_2(T)-\nu_2(\Delta T)\leq cT \big)=
\mathbf{P}\Bigl(\exp\{-r(\nu_2(T)-\nu_2(\Delta T))\}\geq \exp\{-r cT\}\Bigr)
\\
&\leq\frac{\mathbf{E}\exp\{-r(\nu_2(T)-\nu_2(\Delta T))\}}{\exp\{-r cT\}}=
\exp\Bigl\{e^{-r}\frac{\alpha\mu(1-\Delta)}{\lambda+\mu}T-\frac{\alpha\mu(1-\Delta)}{\lambda+\mu}T+r cT\Bigr\}.
\end{aligned}
$$
Choosing $r=-\ln c$ we obtain inequality (\ref{5.4}). $\Box$

\begin{lemma} \label{l5.5} The following inequality holds true
$$\mathbf{P}\bigg(\bigcap\limits_{l=1}^{[cT]}C_l,\overline{B_c},
\sum\limits_{l=1}^{[cT]}\zeta_{k_l}(\xi(\tau_{k_l}-))\leq aT\bigg)\leq
\bigg(\frac{1}{[\delta T]}\bigg)^{[cT]}\exp(aT),
$$
where $C_l$, $1\leq l \leq [cT]$ are defined by (\ref{Cl}) on the previous section.
\end{lemma}

\noindent
\DD. Define random variables
$$
\tilde{\zeta}_{k_l}(m_l):=\left\{
           \begin{array}{cl}
           \zeta_{k_l}(m_l), & \text{ if }\zeta_{k_l}(m_l)\leq [\delta T],\\
           \gamma_{l}, & \text{ if }\zeta_{k_l}(m_l)> [\delta T],
           \end{array}
                           \right.
$$  
where random variables $\gamma_{l}$, $1\leq l \leq [cT]$ are mutually independent and do not depend on
$\zeta_{k_l}(m_l)$, $\xi(\tau_{k_l}-)$, $m_l\in \mathbb{N}$,  $1\leq l \leq [cT]$, $\nu_1(\cdot)$,
$\nu_2(\cdot)$
and $\mathbf{P}(\gamma_{l}=r)=\dfrac{r}{[\delta T]}$, $1\leq r \leq [\delta T]$.

Since $\tilde{\zeta}_{k_l}(\xi(\tau_{k_l}-))\leq \zeta_{k_l}(\xi(\tau_{k_l}-))$, then
\beq \label{5.11}
\mathbf{P}\bigg(\bigcap\limits_{l=1}^{[cT]}C_l,\overline{B_c},
\sum\limits_{l=1}^{[cT]}\zeta_{k_l}(\xi(\tau_{k_l}-))\leq aT\bigg)\leq 
\mathbf{P}\bigg(\bigcap\limits_{l=1}^{[cT]}C_l,\overline{B_c},
\sum\limits_{l=1}^{[cT]}\tilde{\zeta}_{k_l}(\xi(\tau_{k_l}-))\leq aT\bigg).
\eeq
Denote
$$
V:=\Bigl\{v_1\in \mathbb{Z}^+,\dots,v_{[cT]}\in \mathbb{Z}^+:v_1+\dots+v_{[cT]}\leq aT, 
\max\limits_{1 \leq l \leq [cT]}v_l \leq  [\delta T] \Bigr\}.
$$
We have
$$
\begin{aligned}
\mathbf{P} & \bigg(\bigcap\limits_{l=1}^{[cT]}C_l,\overline{B_c},
\sum\limits_{l=1}^{[cT]}\tilde{\zeta}_{k_l}(\xi(\tau_{k_l}-))\leq aT\bigg)
\\
&=\sum\limits_{v_1,\dots,v_{[cT]}\in V}\mathbf{P}\bigg(\bigcap\limits_{l=1}^{[cT]}C_l,
\tilde{\zeta}_{k_1}(\xi(\tau_{k_1}-))=v_1,\dots,\tilde{\zeta}_{k_{[cT]}}(\xi(\tau_{k_{[cT]}}-))=v_{[cT]},\overline{B_c}\bigg)
\\
&=\sum\limits_{v_1,\dots,v_{[cT]}\in V}\mathbf{P}\bigg(\bigcap\limits_{l=1}^{[cT]}C_l,\bigcap\limits_{l=1}^{[cT]}D_l,\overline{B_c}\bigg)
=\sum\limits_{v_1,\dots,v_{[cT]}\in V}\mathbf{P}\bigg(\bigcap\limits_{l=1}^{[cT]}(C_l \cap D_l) \ \bigg| \ \overline{B_c}\bigg)
\mathbf{P}(\overline{B_c}),
\end{aligned}
$$
where
$$
\ D_l:=\big\{\omega:\tilde{\zeta}_{k_l}(\xi(\tau_{k_l}-))=v_l\big\}, \ \ \ 1\leq l \leq [cT].
$$
Let $D_0:=\Omega$, $C_0:=\Omega$.
We will show that for $1 \leq  l \leq [cT]$ the following inequality holds
\beq \label{5.5}
\mathbf{P}\bigg(
D_l,C_l \ \bigg| \bigcap\limits_{d=0}^{l-1}(C_d\cap D_d),\overline{B_c}\bigg)\leq
\frac{1}{[\delta T]}.
\eeq
Note that by definition, the family of random variables
$\tilde{\zeta}_{k_l}(m_l)$, $m_l\in \mathbb{N}$ do not depend on
$\tilde{\zeta}_{k_1}(m_1)$, $m_{1}\in \mathbb{N}$, $\dots$,
  $\tilde{\zeta}_{k_{l-1}}(m_{l-1})$, $m_{l-1}\in \mathbb{N}$
   and  $\xi(\tau_{k_1}-),\dots,\xi(\tau_{k_{l}}-)$, $\nu_2(\cdot)$.
Thus,
$$
\begin{aligned}
& \mathbf{P}\bigg(
D_l,C_l  \ \bigg| \ \bigcap\limits_{d=0}^{l-1}(C_d\cap D_d),\overline{B_c}\bigg)
\\
&=\mathbf{P}\bigg(
D_l  \bigg| \ \bigcap\limits_{d=0}^{l}C_d \bigcap\limits_{d=0}^{l-1}D_d,\overline{B_c}\bigg)
\mathbf{P}\bigg(
C_l  \ \bigg| \ \bigcap\limits_{d=0}^{l-1}(C_d\cap D_d),\overline{B_c}\bigg)
\leq \mathbf{P}\bigg(
D_l  \bigg| \ \bigcap\limits_{d=0}^{l}C_d, \bigcap\limits_{d=0}^{l-1}D_d,\overline{B_c}\bigg)
\\
&=\sum\limits_{r=[\delta T]}^\infty \mathbf{P}\bigg(
\tilde{\zeta}_{k_l}(r)=v_l  \bigg| \   \xi(\tau_{k_l}-)=r, \bigcap\limits_{d=0}^{l}C_d, \bigcap\limits_{d=0}^{l-1}D_d,\overline{B_c}\bigg)
 \mathbf{P}\bigg(\xi(\tau_{k_l}-)=r  \ \bigg| \ \bigcap\limits_{d=0}^{l}C_d, \bigcap\limits_{d=0}^{l-1}D_d,\overline{B_c}\bigg)
\\
&=\sum\limits_{r=[\delta T]}^\infty 
\frac{1}{[\delta T]} \mathbf{P}\bigg(\xi(\tau_{k_l}-)=r  \ \bigg| \ \bigcap\limits_{d=0}^{l}C_d, \bigcap\limits_{d=0}^{l-1}D_d,\overline{B_c}\bigg)=\frac{1}{[\delta T]},
\end{aligned}
$$
where we used the fact that, if $r\geq [\delta T]$, then  
$\mathbf{P}(\tilde{\zeta}_{k_l}(r)=v_l)=\dfrac{1}{[\delta T]}$.
Using (\ref{5.5}) we obtain
$$
\mathbf{P}\bigg(\bigcap\limits_{l=1}^{[cT]}C_l,\bigcap\limits_{l=1}^{[cT]}D_l,\overline{B_c}\bigg)=
\prod\limits_{l=1}^{[cT]}\mathbf{P}\bigg(
D_l,C_l \ \bigg| \bigcap\limits_{d=0}^{l-1}(C_d\cap D_d),\overline{B_c}\bigg)
\mathbf{P}(\overline{B_c})\leq
\bigg(\dfrac{1}{[\delta T]}\bigg)^{[cT]}.
$$
Thus, 
$$
\mathbf{P}\bigg(\bigcap\limits_{l=1}^{[cT]}C_l,\overline{B_c},
\sum\limits_{l=1}^{[cT]}\tilde{\zeta}_{k_l}(\xi(\tau_{k_l}-))\leq aT\bigg)\leq
\sum\limits_{v_1,\dots,v_{[cT]}\in V}\bigg(\dfrac{1}{[\delta T]}\bigg)^{[cT]}=
\mathbf{P}\bigg(
\sum\limits_{l=1}^{[cT]}\zeta_{k_l}([\delta T])\leq aT\bigg).
$$
Therefore, from (\ref{5.11}) and Lemma~\ref{l5.2} it follows that
$$\mathbf{P}\bigg(\bigcap\limits_{l=1}^{[cT]}C_l,\overline{B_c}
\sum\limits_{l=1}^{[cT]}\zeta_{k_l}(\xi(\tau_{k_l}-))\leq aT\bigg)
\leq\bigg(\frac{1}{[\delta T]}\bigg)^{[cT]}\exp(aT).
$$
$\Box$

\section*{Acknowledgments}
This work is supported by FAPESP grant 2017/20482-0.

AL supported by RSF according to the research project 18-11-00129. AL thanks Institute of Mathematics and 
Statistics of University of S\~ao Paulo for hospitality. AY thanks CNPq and  
FAPESP for the financial support via grants 301050/2016-3 and 2017/10555-0, respectively.

\end{document}